\documentclass{article}

\usepackage{arxiv}

\usepackage[utf8]{inputenc} 
\usepackage[T1]{fontenc}    
\usepackage{hyperref}       
\usepackage{url}            
\usepackage{booktabs}       
\usepackage{amsfonts}       
\usepackage{nicefrac}       
\usepackage{microtype}      
\usepackage{graphicx}
\usepackage{doi}
\usepackage{amsthm}
\usepackage{dsfont}
\usepackage{tikz}
\usepackage{mathtools}
\usepackage{tcolorbox}
\usepackage{amsmath}
\usepackage{nicefrac}
\usepackage{physics}
\usepackage{bbm}
\usepackage{wasysym}

\newtheorem{theorem}{Theorem}
\newtheorem{lemma}{Lemma}

\newtheorem{example}{Example}
\theoremstyle{definition}

\DeclareMathOperator{\sgn}{sgn}
\DeclareMathOperator{\pvv}{p.\hspace{-0.1cm}v.}

\renewcommand{\phi}{\varphi}
\newcommand{\mi}{\mathrm{i}}
\renewcommand{\H}{\mathcal{H}}

\newcommand{\R}{\mathbb{R}}
\newcommand{\Z}{\mathbb{Z}}
\newcommand{\N}{\mathbb{N}}

\title{Some properties of a modified Hilbert transform}

\author{\href{https://orcid.org/0000-0002-2577-1421}{\includegraphics[scale=0.06]{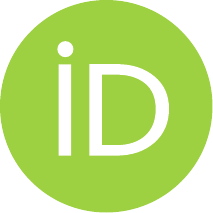}\hspace{1mm}Matteo ~Ferrari} \\
	Faculty of Mathematics \\
    University of Vienna\\
    1090 Vienna, Austria \\
	\texttt{matteo.ferrari@univie.ac.at}}

\renewcommand{\headeright}{}
\renewcommand{\undertitle}{}
\renewcommand{\shorttitle}{}

\date{}

\hypersetup{
pdftitle={hilbert},
pdfsubject={q-bio.NC, q-bio.QM},
pdfauthor={Ferrari M.~Kubin A.},
pdfkeywords={Hilbert transform},
}

\begin{document}
\maketitle

\begin{abstract}
Recently, Steinbach et al. introduced a novel operator $\H_T: L^2(0,T) \to L^2(0,T)$, known as the modified Hilbert transform. This operator has shown its significance in space-time formulations related to the heat and wave equations. In this paper, we establish a direct connection between the modified Hilbert transform $\H_T$ and the canonical Hilbert transform $\H$. Specifically, we prove the relationship $\H_T \phi = -\H \widetilde{\phi}$, where $\phi \in L^2(0,T)$ and $\widetilde{\phi}$ is a suitable extension of $\phi$ over the entire $\R$. By leveraging this crucial result, we derive some properties of $\H_T$, including a new inversion formula, that emerge as immediate consequences of well-established findings on $\H$.
\end{abstract}

\section{Introduction and main result}
In \cite{SteinbachZank2020}, a modified Hilbert transform $\mathcal{H}_T$, associated with the bounded interval $(0,T)$, has been defined. Given the Fourier series of $\phi \in L^2(0,T)$
\begin{equation*}
	\phi(t) = \sum_{k=0}^\infty \phi_k \sin \left(\left( \frac{\pi}{2} + k \pi\right)\frac{t}{T} \right) \quad \text{with} \quad \phi_k =  \frac{2}{T} \int_0^T \phi(s) \sin \left( \left( \frac{\pi}{2} + k \pi\right) \frac{s}{T} \right) \dd s,
\end{equation*}
the modified Hilbert transform is defined as
\begin{equation} \label{defHilbZank}
	\H_T \phi(t) = \sum_{k=0}^\infty \phi_k \cos \left(\left( \frac{\pi}{2} + k \pi\right)\frac{t}{T} \right) \quad t \in (0,T).
\end{equation}

This operator has been employed in the context of space-time discretizations of PDEs using both finite element \cite{HauserZank2024,PerugiaSchwabZank2023,SteinbachZank2020} and boundary element \cite{SteinbachUrzuaTorresZank2022} methods. It is particularly well-suited because, by defining the Sobolev space 
\begin{equation*}
	H_{0,}^{\nicefrac{1}{2}}(0,T) = \{ \phi \in H^{\nicefrac{1}{2}}(0,T) : \phi(0)=0\},
\end{equation*}
it can be shown (see \cite{SteinbachZank2020}) that $-\partial_t \H_T : H_{0,}^{\nicefrac{1}{2}}(0,T) \to [H_{0,}^{\nicefrac{1}{2}}(0,T)]'$ induces an equivalent norm in $H_{0,}^{\nicefrac{1}{2}}(0,T)$. Furthermore, it holds $\langle \phi , \H_T \phi \rangle_{L^2(0,T)} \ge 0$ for all $\phi \in L^2(0,T)$.

An alternative integral representation of $\H_T$ suited for numerical schemes has been presented in \cite[Lemma 2.1]{SteinbachZank2021}.
\begin{lemma}
For $\phi \in L^2(0,T)$, the operator $\H_T$ allows the integral representation
\begin{equation} \label{HT}
	\H_T \phi(t) = \frac{1}{2T} \pvv \int_0^T \phi(s) \left[ \csc \left( \frac{\pi(s+t)}{2T} \right) + \csc \left( \frac{\pi(s-t)}{2T} \right)\right] \dd s, \quad t \in (0,T)
\end{equation}
as Cauchy principal value integral.
Moreover, if $v \in H^1(0,T)$ it holds
\begin{equation} \label{SZ}
	\begin{aligned}
	\H_T \phi(t) = & -\frac{2}{\pi} \phi(0)  \log \left( \tan \left(\frac{\pi t}{4T}\right)\right) - \frac{1}{\pi} \int_0^T \partial_t \phi(s) \log\left( \tan \left( \frac{\pi(s+t)}{4T} \right) \tan \left( \frac{\pi|s-t|}{4T} \right)\right) \dd s,
	\end{aligned}
\end{equation}
for $t \in (0,T)$ as a weakly singular integral.
\end{lemma}
The Hilbert transform $\H$ of a function $\phi$ is defined as Cauchy principal value integral
\begin{equation*}
	\H \phi(t) = \frac{1}{\pi} \pvv\int_{\R} \frac{\phi(s)}{t-s} \, ds,
\end{equation*}
whenever it exists (see \cite[Chapter 9]{ButzerNessel1971}, \cite{King2009} and references therein).

\begin{figure}[h!]
\centering
  \includegraphics[width=0.8\linewidth]{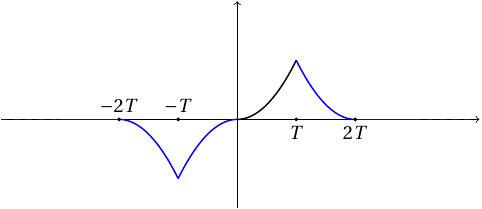}
\caption{For $\phi(x) = x^2$, $\phi : [0,T] \to \R$, it is plotted $\bar \phi : \R \to \R$ defined in \eqref{barphi}.}
\label{fig1}
\end{figure}

The following relationship between the Hilbert transform $\H$ and its modified version $\H_T$ has been established in \cite[Theorem 4.3]{SteinbachMissoni2023}.
\begin{theorem}
For $\phi \in L^2(0,T)$, it holds
\begin{equation*}
	\H_T \phi = -\H \bar \phi + B \phi \quad \quad \text{in~} L^2(0,T),
\end{equation*}
where $B : L^2(0,T) \to L^2(0,T)$ is a compact operator, and $\bar \phi$ is the reflection (see Figure \ref{fig1})
\begin{equation} \label{barphi}
	\bar \phi(s) =
	\begin{cases}
			-\phi(s+2T) & s \in (-2T,-T), \\
					-\phi(-s) & s \in (-T,0), \\
		\phi(s) & s \in (0,T), \\
		\phi(2T-s) & s \in (T,2T), \\
		0 & \text{else}.
	\end{cases}
\end{equation}
\end{theorem}

In this paper, we prove that $\H_T$ is, in fact, precisely the Hilbert transform $\H$ applied to a specific odd periodic extension with alternating signs.
\begin{theorem} \label{mainth}
For $\phi \in L^2(0,T)$, it holds
\begin{equation*}
	\H_T \phi = -\H \widetilde \phi \quad \quad \text{in~} L^2(0,T),
\end{equation*}
where $\widetilde{\phi}$ is the periodic reflection (see Figure \ref{fig2})
\begin{equation} \label{tphi}
	\widetilde \phi(s) =
	\begin{cases}
		-\phi(s+2T-4kT) & s \in ((4k-2)T,(4k-1)T), \, \,  k \in \Z, \\
		-\phi(4kT-s) & s \in ((4k-1)T,4kT), \, \,  k \in \Z, \\
		\phi(s-4kT) & s \in (4kT,(4k+1)T), \, \, k \in \Z,  \\
		\phi(4kT+2T-s) & s \in ((4k+1)T,(4k+2)T), \, \, k \in \Z.  \\
	\end{cases}
\end{equation}
\begin{figure}[h!]
\centering
  \includegraphics[width=0.8\linewidth]{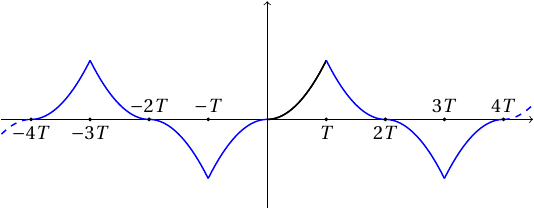}
\caption{For $\phi(x) = x^2$, $\phi : [0,T] \to \R$, it is plotted $\widetilde \phi : \R \to \R$ defined in \eqref{tphi}.}
\label{fig2}
\end{figure}
\end{theorem}
\begin{example}
From \eqref{SZ}, we obtain that for $\phi(s) = 1$ with $s \in (0,T)$, we have $\H_T \phi(t) = -\frac{2}{\pi} \log\left| \tan \left( \frac{\pi t}{4 T} \right) \right|$. It can be readily seen that, for $s \in \R$, $\widetilde{\phi}(s) = \sgn(\sin(a s))$ with $a = \frac{\pi}{2T}$, and for this function, it is well-known that $\H \widetilde{\phi} (t) = \frac{2}{\pi} \log\left| \tan \left(\frac{at}{2}\right) \right|$ (see \cite[Equation 6.9]{King2009}).
\end{example}
\begin{example}
A function that satisfies a periodicity of the form \eqref{tphi} is $\phi(s) = \sin(\pi s)$ with $T = \frac{1}{2}$. In this case, on $\mathbb{R}$, we have the trivial extension $\widetilde{\phi}(s) = \sin(\pi s)$. It is well-known that $\H \widetilde{\phi}(t) = -\cos(\pi t)$ (see \cite[Equation 1.8]{King2009}). From \eqref{HT}, we can calculate
\begin{align*}
	\H_T \phi(t) & = \pvv \int_0^{\frac{1}{2}} \sin(\pi s) \left[ \csc\left( \pi(s+t) \right) + \csc\left( \pi(s-t) \right) \right] \dd s 
	\\ &= \pvv \int_t^{t+\frac{1}{2}} \frac{\sin(\pi(s-t))}{\sin(\pi s)} \dd s + \pvv \int_{-t}^{\frac{1}{2}-t} \frac{\sin(\pi(s-t))}{\sin(\pi s)} \dd s = \cos(\pi t).
\end{align*}
\end{example}

\section{Proofs of the main result}
In this section, we present three distinct proofs of Theorem \ref{mainth}. The first proof relies on the integral representation \ref{HT}, while the second proof utilizes an alternative definition of the Hilbert transform $\mathcal{H}$ specifically designed for periodic functions. Finally, the third proof employs Fourier series.
\subsection{Proof based on the integral representation}
The first proof is established by leveraging the Laurent series expansion of the cosecant function, as presented in \cite[Formula (9.3.30)]{Hille1959}
\begin{equation} \label{MitLeg}
	\csc(z) = \frac{1}{z} + 2z \sum_{k = 1}^{\infty} \frac{(-1)^k}{z^2 - k^2 \pi^2}.
\end{equation}
This series converges absolutely and uniformly for $|z| < \pi$, in fact for each $N \in \N$ it holds
\begin{equation*} 
	\left| \csc(z) - \frac{1}{z} \right| \le \frac{2|z|}{\pi^2} \sum_{k = 1}^{N} \frac{1}{\left| \left( \frac{z}{k\pi}\right)^2 - 1\right|} \frac{1}{k^2} \le \frac{|z|}{3}.
\end{equation*}
Starting from the definitions of Hilbert transform $\H$ and $\widetilde \phi$ in \eqref{tphi}, we can write
\begin{align*}
	 \H \widetilde \phi(t) &  = \frac{1}{\pi} \left[ \ldots + \pvv \int_{-3T}^{-T} \frac{\widetilde \phi(s)}{t-s} \dd s + \pvv \int_{-T}^{T} \frac{\widetilde \phi(s)}{t-s} \dd s + \pvv \int_{T}^{3T} \frac{\widetilde \phi(s)}{t-s} \dd s  +\ldots \right] 
	\\ & = \frac{1}{\pi} \left[ \ldots - \pvv \int_{-T}^{T} \frac{\widetilde \phi(s)}{t-s-2T} \dd s + \pvv \int_{-T}^{T} \frac{ \widetilde \phi(s)}{t-s} \dd s - \pvv \int_{-T}^{T} \frac{\widetilde  \phi(s)}{t-s+2T} \dd s  +\ldots \right]
\end{align*}
where we have used the alternating signs periodicity of $\widetilde{\phi}$. Therefore, we can compactly write
\begin{align*}
	\H \widetilde \phi(t) & = \frac{1}{\pi} \sum_{k=-\infty}^{+\infty} \pvv \int_{-T}^T \widetilde{\phi}(s) \frac{(-1)^k }{t-s+2kT} \dd s.
\end{align*}
By interchanging the integral and the summation (which is possible due to the uniform convergence of the series), and utilizing \eqref{MitLeg}, we obtain
\begin{align*}
	\H \widetilde \phi(t) & = \frac{1}{\pi} \pvv \int_{-T}^T \widetilde \phi(s) \left[ \sum_{k=-\infty}^{+\infty} \frac{(-1)^k}{t-s+2kT} \right] \dd s
\\ & = \frac{1}{\pi} \pvv \int_{-T}^T \widetilde \phi(s) \left[\frac{1}{t-s} + 2(t-s) \sum_{k=1}^{+\infty} \frac{(-1)^k}{(t-s)^2-4k^2T^2} \right] \dd s
	\\ & = \frac{1}{2T} \pvv \int_{-T}^T \widetilde \phi(s) \csc\left(\frac{\pi(t-s)}{2T}\right) \dd s.
\end{align*}
From the latter, we obtain the modified Hilbert transform $\H_T$ in integral form \eqref{HT} by recalling the definition \eqref{tphi}, and writing
\begin{align*}
	\H \widetilde \phi(t) & = - \frac{1}{2T} \pvv \int_{-T}^0 \phi(-s) \csc\left(\frac{\pi(t-s)}{2T}\right) \dd s + \frac{1}{2T} \pvv \int_0^T \phi(s) \csc\left(\frac{\pi(t-s)}{2T}\right) \dd s
	\\ & = - \frac{1}{2T} \pvv \int_0^T \phi(s) \csc\left(\frac{\pi(t+s)}{2T}\right) \dd s - \frac{1}{2T} \pvv \int_0^T \phi(s) \csc\left(\frac{\pi(t-s)}{2T}\right) \dd s
	\\ & = -\H_T \phi(t).
\end{align*}
\subsection{Proof based on the Hilbert transform for periodic functions}
Given $\phi \in L^2(0,T)$, we observe that $\widetilde{\phi}$ as in \eqref{tphi} is actually periodic with a period of $4T$, and it is known that the Hilbert transform $\H$ for a periodic function can be calculated using (see \cite{Pandey1997}):
\begin{equation} \label{eqper}
	\H \widetilde \phi (t) = \frac{1}{4T} \pvv \int_{-2T}^{2T} \widetilde \phi(s) \cot\left(\frac{\pi(t-s)}{4T} \right) \dd s.
\end{equation}
Let us expand the four contributions of the integral above to recover $\H_T \phi$
\begin{align*}
	\H \widetilde \phi (t) & = -\frac{1}{4T} \left[ \pvv \int_{-2T}^{-T} \phi(s+2T) \cot\left(\frac{\pi(t-s)}{4T} \right) \dd s - \pvv \int_{-T}^{0}  \phi(-s) \cot\left(\frac{\pi(t-s)}{4T} \right) \dd s \right.
	\\ & \hspace{0.45cm} \left. + \pvv \int_{0}^{T} \phi(s) \cot\left(\frac{\pi(t-s)}{4T} \right) \dd s + \pvv \int_{T}^{2T} \phi(2T-s) \cot\left(\frac{\pi(t-s)}{4T} \right) \dd s \right]
	\\ & = -\frac{1}{4T} \pvv \int_{0}^{T} \phi(s) \cot\left(\frac{\pi (t-s+2T)}{4T} \right) \dd s - \frac{1}{4T} \pvv \int_{0}^{T}  \phi(s) \cot\left(\frac{\pi(t+s)}{4T} \right) \dd s
	\\ & \hspace{0.45cm} + \frac{1}{4T} \pvv \int_{0}^{T} \phi(s) \cot\left(\frac{\pi(t-s)}{4T} \right) \dd s + \frac{1}{4T} \pvv \int_{0}^{T} \phi(s) \cot\left(\frac{\pi(t+s-2T)}{4T} \right) \dd s.
\end{align*}
Finally, using the trigonometric formula 
\begin{align} \label{trigon}
	2\csc (2x) =  -\cot\left(x \pm \frac{\pi}{2} \right)  + \cot(x), \quad |x| < \frac{\pi}{2},
\end{align}
we again conclude that $ \H_T \phi = - \H \widetilde \phi$ in $L^2(0,T)$.

\subsection{Proof based on Fourier series}

For $f$ with period $2T$ in $L^2(-T,T)$, the circular Hilbert transform is defined as (see \cite[Section 6.4]{King2009})
\begin{equation*}
	\H f(t) = \mi \sum_{k=1}^{\infty} \left[ f_{-k} e^{-\mi k \pi \frac{t}{T}} - f_k e^{\mi k \pi \frac{t}{T}}\right] \quad \text{with} \quad 	f_k = \frac{1}{2T} \int_{-T}^T f(s) e^{-\mi k \pi \frac{s}{T}} \dd s.
\end{equation*}
Let us verify that this definition is consistent with the original definition of the modified Hilbert transform $\H_T$ by Steinbach and Zank \eqref{defHilbZank}. That is, we will demonstrate that for a given function $\phi \in L^2(0,T)$, we have $\mathcal{H} \widetilde{\phi} = -\mathcal{H}_T \phi$ in $L^2(0,T)$, with $\widetilde \phi$ as in \eqref{tphi}. Recalling that $\widetilde \phi$ is a periodic function with period $4T$, let us begin by writing
\begin{equation*}
	\H \widetilde \phi(t) = \mi \sum_{k=1}^{\infty} \left[ \widetilde{\phi}_{-k} e^{-\mi k \pi \frac{t}{2T}} - \widetilde{\phi}_k e^{\mi k \pi \frac{t}{2T}}\right],
\end{equation*}
and readily one can verify
\begin{align*}
	\widetilde{\phi}_k = \frac{1}{4T} \int_{-2T}^{2T} \widetilde \phi(s) e^{-\mi k \pi \frac{s}{2T}} \dd s & = \frac{1}{4T} \int_0^{2T} \widetilde\phi(s) \left[e^{-\mi k \pi \frac{s}{2T}} - e^{\mi k \pi \frac{s}{2T}}\right] \dd s
	\\ & = - \frac{\mi}{2T} \int_0^{2T} \widetilde\phi(s) \sin\left(k \pi \frac{s}{2T}\right) \dd s = -\widetilde{\phi}_{-k}.
\end{align*}
Hence, we continue
\begin{align*}
	\H \widetilde \phi(t) & = - \sum_{k=1}^{\infty} \left[e^{-\mi k \pi \frac{t}{2T}} + e^{\mi k \pi \frac{t}{2T}}\right] \frac{1}{2T} \int_0^{2T} \widetilde \phi(s) \sin\left(k \pi \frac{s}{2T}\right) \dd s
	\\ & = -\sum_{k=1}^{\infty} 2 \cos\left(k \pi \frac{t}{2T} \right) \frac{1}{2T} \int_0^{2T} \widetilde\phi(s) \sin\left(k \pi \frac{s}{2T}\right) \dd s.
\end{align*}
Moreover, we can now write
\begin{align*}
	\frac{1}{2T} \int_0^{2T} \widetilde\phi(s) \sin\left(k \pi \frac{s}{2T}\right) \dd s  & = \frac{1}{2T} \int_0^T \phi(s) \sin\left(k \pi \frac{s}{2T}\right) \dd s+\frac{1}{2T} \int_T^{2T} \phi(2T-s) \sin\left(k \pi \frac{s}{2T}\right) \dd s 
	\\ & = \frac{1}{2T} \int_0^T \phi(s) \left[\sin\left(k \pi \frac{s}{2T}\right) - \sin(k \pi) \sin\left(k \pi \frac{s}{2T}\right)\right] \dd s 
	\\ & = 
	\begin{cases}
		0 & k \text{~ even}, \\
		\displaystyle{\frac{1}{T} \int_0^T \phi(s) \sin\left(k \pi \frac{s}{2T}\right) \dd s} & k \text{~ odd}.
	\end{cases}
\end{align*}
Finally, we conclude
\begin{equation*}
	\H \widetilde \phi(t) = -\sum_{k=0}^{\infty} \cos\left((2k+1) \pi \frac{t}{2T}\right) \frac{2}{T} \int_0^T \widetilde\phi(s) \sin\left((2k+1) \pi \frac{s}{2T}\right) \dd s = -\H_T \phi(t).
\end{equation*}

\section{Consequences of the main result}
In this section we show simple consequences of Theorem \ref{mainth}.

\subsection{Inversion formula}
For $f \in L^2(-T,T)$ and periodic with period $2T$, the inversion formula holds (see \cite[Formula (6.35)]{King2009})
\begin{equation} \label{invH}
	\H^2 f(t) = -f(t) + \frac{1}{2T} \int_{-T}^T f(s) \dd s, \quad \quad \text{in~} L^2(-T,T).
\end{equation}
For $\phi \in L^2(0,T)$, we can calculate
\begin{equation*} 
	\H (\H_T \phi)(t) = -\H^2 \widetilde{\phi}(t) = \widetilde \phi (t) - \frac{1}{4T} \int_{-2T}^{2T} \widetilde \phi(s) \dd s = \phi(t), \quad \quad \text{in~} L^2(0,T),
\end{equation*}
since $\widetilde \phi$ is an odd function.

\subsection{Alternative formula}
For an odd function $f \in L^2(-T,T)$ with period $2T$, it can be shown that
\begin{equation*}
	\H f(t)= \frac{1}{T} \pvv \int_0^T f(s) \frac{\sin\left(\pi \frac{s}{T} \right)}{\cos\left( \pi \frac{s}{T}\right)-\cos\left(\pi \frac{t}{T} \right)} \dd s, \quad \quad \text{in~} L^2(-T,T).
 \end{equation*}
Therefore, we obtain the alternative formula in $L^2(0,T)$
\begin{align*}
	\H_T \phi(t)  = - \H \widetilde \phi(t) & = - \frac{1}{2T} \pvv \int_0^{2T} \widetilde \phi(s) \frac{\sin\left( \pi \frac{s}{2T} \right)}{\cos\left( \pi \frac{s}{2T}\right)-\cos\left( \pi \frac{t}{2T} \right)} \dd s
	\\ & = -\frac{\cos\left( \pi \frac{t}{2T} \right)}{T} \pvv \int_0^T  \phi(s) \frac{\sin\left( \pi \frac{s}{2T} \right)}{\cos^2\left( \pi \frac{s}{2T} \right)-\cos^2\left( \pi \frac{t}{2T} \right)} \dd s.
 \end{align*}
This formula can also be deduced from equation \eqref{HT} using trigonometric identities.

\subsection{Integral representation}
Let suppose that $\phi \in H^1(0,T)$, then the derivative of $\widetilde \phi$ is 
\begin{equation} \label{dphit}
	\partial_t \widetilde \phi(s)_{|_{(-2T,2T)}} =
	\begin{cases}
		-\partial_t \phi(s+2T) & s \in (-2T,-T), \\
		\partial_t \phi(-s) & s \in (-T,0), \\
		\partial_t \phi(s) & s \in (0,T), \\
		-\partial_t \phi(2T-s) & s \in (T,2T).
	\end{cases}
\end{equation}
We note that if $\phi \in H^1(0,T)$, the only possible discontinuity of $\widetilde \phi_{|_{(-2T,2T)}}$ is in the point $s=0$.

Starting from \eqref{eqper}, integrating by parts in the intervals of continuity $(-2T,0)$ and $(0,2T)$ of $\widetilde \phi$, we obtain
\begin{align*}
	 \H_T \phi(t) = -\H \widetilde \phi (t) & = -\frac{1}{\pi} \int_{-2T}^{2T} \partial_t \widetilde \phi(s) \log \left( \sin\left(\frac{\pi|t-s|}{4T} \right) \right) \dd s 
	\\ & \hspace{0.35cm} + \frac{1}{\pi} \left[ (\widetilde \phi(0^-)-\widetilde \phi(0^+)) \log \left( \sin \left( \frac{\pi t}{4T} \right) \right) + (\widetilde \phi(2T^-)-\widetilde \phi(-2T^+)) \log \left( \cos \left( \frac{\pi t}{4T} \right) \right) \right]
	\\ & = -\frac{1}{\pi}  \int_{-2T}^{2T} \partial_t \widetilde \phi(s) \log \left( \sin\left(\frac{\pi|t-s|}{4T} \right) \right) \dd s - \frac{2}{\pi} \phi(0) \log \left( \tan \left( \frac{\pi t}{4T} \right) \right).
\end{align*}
This formula, in the case $\phi(0)=0$ is the original form in which Hilbert wrote the Hilbert transform in \cite{Hilbert1904}. Let split the integral in the four usual intervals, and use \eqref{dphit},
\begin{align*}
	 \H_T \phi (t) & =  \frac{1}{\pi} \int_0^T \partial_t \phi(s) \log \left( \cos\left(\frac{\pi(t+s)}{4T} \right) \right) \dd s - \frac{1}{\pi} \int_0^T  \partial_t \phi(s) \log \left( \sin\left(\frac{\pi(t+s)}{4T} \right) \right) \dd s
	\\ & \hspace{0.3cm} - \frac{1}{\pi} \int_0^T  \partial_t \phi(s) \log \left( \sin\left(\frac{\pi|t-s|}{4T} \right) \right) \dd s + \frac{1}{\pi} \int_{0}^{T} \partial_t \phi(s) \log \left( \cos\left(\frac{\pi(t-s)}{4T} \right) \right) \dd s
	\\ & \hspace{0.3cm} -  \frac{2}{\pi} \phi(0) \log \left( \tan \left( \frac{\pi t}{4T} \right) \right)
	\\ & = - \frac{1}{\pi} \int_0^T \partial_t \phi(s) \log \left( \tan\left(\frac{\pi(t+s)}{4T}\right) \tan \left( \frac{\pi| t-s|}{4T} \right) \right) \dd s -  \frac{2}{\pi} \phi(0) \log \left( \tan \left( \frac{\pi t}{4T} \right) \right).
\end{align*}
We have obtained exactly the result of Lemma \ref{SZ} in an alternative way.
\section{Acknowledgements}
This work is partially supported by the INdAM-GNCS project ``Metodi numerici per lo studio di strutture geometriche parametriche complesse'' (CUP E53C22001930001) and by the MIUR project ``Dipartimenti di Eccellenza 2018-2022'' (CUP E11G18000350001). Moreover, this research was funded in part by the Austrian Science Fund (FWF) project \href{https://doi.org/10.55776/F65}{10.55776/F65}

\bibliographystyle{plain}
\bibliography{bibtex_num}

\end{document}